\documentclass[11pt]{amsart}
\usepackage{amssymb}

\def \R {\mathbb R}

\def \gg {\mathfrak g}

\def \gu {\mathfrak u}
\def \gv {\mathfrak v}

\def \gI {\mathfrak I}
\def \gJ {\mathfrak J}

\def \gW {\mathfrak W}

\def \CC {\mathcal C}

\def \CF {\mathcal F}

\def \CL {\mathcal L}

\def \CS {\mathcal S}

\def \CU {\mathcal U}

\def \bC {\mathbf C}
\def \bD {\mathbf D}

\def \bd {\mathbf d}

\def \h {\hbar}

\def \norm {\rm {norm}}

\def \inv {^{-1}}
\def \o {\overline}
\def \wt {\widetilde}

\def \ad {\mathop {\rm ad}{}}

\def \adj {\mathop {\rm adj}{}}
\def \tr {\mathop {\rm tr}}

\def \supp {\mathop {\rm supp}}
\def \norm#1{\Vert #1 \Vert}
\def \sl {\mathfrak {sl}}

\newtheorem{theorem}{Theorem}[section]

\newtheorem{lemma}[theorem]{Lemma}
\newtheorem{proposition}[theorem]{Proposition}

\newtheorem{proposition-definition}[theorem]{Proposition and
Definition}

\begin{document}

\title[invariant distributions]{Convolution of Invariant
Distributions:\\Proof of the Kashiwara-Vergne conjecture}

\author{Martin Andler}
\address{D\'{e}partement de Math\'ematiques,
Universit\'{e} de Versailles Saint Quentin, 78035
Versailles C\'{e}dex}
\email{andler@math.uvsq.fr}
\thanks{The research of the first author was supported by CNRS (EP1755)}
\author{Siddhartha Sahi}
\address{Mathematics Department, Rutgers
University, New Brunswick, NJ 08903}
\email{sahi@math.rutgers.edu}
\thanks{The research of the second author was supported by an NSF grant}

\author {Charles Torossian}
\address {D\'epartement de math\'ematiques et applications,
\'Ecole normale sup\'erieure, 45 rue d'Ulm, 75230 Paris CEDEX 05}
\email{Charles.Torossian@ens.fr}
\thanks{The research of the third author was supported by CNRS (UMR 8553)}

\begin{abstract} Consider the Kontsevich $\star$-product on the symmetric
algebra of a finite dimensional Lie algebra $\gg$, regarded as the
algebra of distributions with support $0$ on $\gg$. In this paper, we
extend this
$\star$-product to distributions satisfying an appropriate 
support condition. As a consequence, we prove a long
standing conjecture of Kashiwara-Vergne on the convolution of germs
of invariant distributions on the Lie group $G$.
\end{abstract}

\maketitle

\section{Introduction}
In several problems in harmonic analysis on Lie groups, one
needs to relate invariant distributions on a Lie
group $G$ to invariant 
distributions on its Lie algebra $\gg$. For instance it is a
central aspect in Harish Chandra's work in the semi-simple
case. The symmetric algebra $\CS(\gg)$ and the enveloping algebra 
$\CU(\gg)$ can be regarded as convolution algebras of distributions supported at
$O$ in $\gg$ and $1$ in $G$ respectively. Therefore the
Harish Chandra isomorphism between the ring of invariants in
$\CS(\gg)$ and
$\CU(\gg)$ can be seen in
this light. 
At a more profound level, Harish Chandra's regularity result for
invariant eigendistributions on the group involves the lifting of the corresponding
result on the Lie algebra.  
\par
Using the orbit method of Kirillov, Duflo
defined an isomorphism extending the Harish
Chandra homomorphism to the case of general Lie groups
(\cite{duflo-70},
\cite{duflo-71}, \cite{duflo-77}). This result was crucial
in the proof of local solvability of invariant differential
operators on Lie groups, established by Ra\"\i s 
\cite{rais-71} for nilpotent groups, then by Duflo and
Ra\"\i s 
\cite{duflo-rais} for solvable groups, and finally by Duflo
\cite{duflo-77} in the general case.
\par Soon thereafter
Kashiwara and Vergne \cite{kash-ver} conjectured that a
natural extension of the Duflo isomorphism to
germs of distributions on $\gg$, when restricted to
invariant germs with appropriate support, should carry the
convolution on
$\gg$ to the convolution on $G$. They proved this conjecture for
solvable groups, and showed that in general their
conjecture implies the local solvability result mentioned
above. (See also an observation of Ra\"\i s in
\cite{rais-review}.) 
\par
Apart from the special case of $\sl(2,\R)$ considered by Rouvi\`ere 
\cite{rouv-81}, the Kashiwara-Vergne
conjecture resisted all attempts 
until 1999, when it was established for arbitrary groups, but
under the restriction that one of the distributions have
point support \cite{ads, ads-cras}. This suffices for
many applications, including the local solvability result mentioned
above. Around the same time, Vergne \cite{vergne} proved the conjecture
for arbitrary germs, but for a special class of Lie algebras, the {\em
quadratic} Lie algebras (those admitting an invariant
non-degenerate quadratic form).
\par
In this paper, we prove the
Kashiwara-Vergne conjecture in
full generality, using, as in
\cite{ads}, the Kontsevich quantization of
the dual of a Lie algebra.
\par
Let us now outline the result and our method. Let $G$ be a Lie
group with Lie algebra $\gg$, and $\exp : \gg \to G$ the
exponential map. Let $q$ be the analytic function on $\gg$ defined by
\begin{equation} q(X) = \det \left(\frac{e^{\ad X/2} -
e^{-\ad X/2}}{\ad X}\right) ^{1/2}. \end{equation}
Define, for a distribution $u$ on $\gg$
\begin{equation} \eta (u) = \exp_*(u . q), \end{equation}
where $\exp_*$ is the pushforward of distributions under the exponential
map. We consider germs at $0$ of distributions on $\gg$, henceforth
simply {\em germs}. For distributions $u,v, \dots$ the corresponding
germs will be denoted $\mathfrak {u,v},\dots$; we use the same
notation $\eta$ for the induced map on germs. 

We wish to consider the convolution of germs. This notion is well
defined under a certain asymptotic support condition which we now
describe. If $U$ is a subset of $\R^n$, one defines the asymptotic cone
of $U$ at $x \in \R^n$ as the set $C_x(U)$ of limit points of all
sequences 
\begin{equation} a_n(x_n -x) \end{equation}
for all sequences $x_n \to x$ and all sequences $a_n
\in \R^+$. Clearly, $C_x(U)$ depends only on the intersction of $U$
with an arbitrarily small neighborhood of $x$. If
$M$ is a manifold, and
$x \in M$, $U \subset M$, using a coordinate chart, one can once again
define $C_x(U)$ as a cone in the tangent space $T_xM$. 

Let $u$ be a distribution on $\gg$. Then $C_0(\supp u)$ depends only
on the germ $\gu$, and we will write it as $C_0[\gu]$.  Assume
that
$\gu,\gv$ are germs such that
\begin{equation} \label{=KV-condition} C_0[\gu] \cap -
C_0[\gv] =\{0\}.
\end{equation}
When two germs verify (\ref{=KV-condition})
we will say that they are {\em compatible}. In this case the
(abelian) convolution on the Lie algebra  $\gu\ast_{\gg} \gv$ is well
defined as a germ on $\gg$. Also, using the
Campbell-Hausdorff formula, it is easy to see that the Lie group 
convolution $\eta(\gu)\ast_{G} \eta(\gv)$ is well defined as the
germ at $1$ of a distribution on the group $G$. 

The Lie algebra $\gg$ acts on functions on $\gg$ by adjoint vector
fields. The dual action descends to germs. We call a germ {\em
invariant} if it is annihilated by all elements of the Lie algebra.
Our main theorem is~:

\begin{theorem} [Kashiwara-Vergne conjecture] Assume that
$\gu$ and $\gv$ are compatible invariant germs. Then
\begin{equation} \eta(\gu \ast_{\gg} \gv) = \eta(\gu) \ast_G
\eta(\gv).
\end{equation}
\label{main-theorem}
\end{theorem}
As in
\cite{kash-ver}, one can reformulate the conjecture slightly by
considering, for $t \in \R$, the Lie group
$G_t$ with Lie algebra $\gg_t$, where $\gg_t$ is $\gg$ as a vector
space, equipped with the Lie bracket 
\begin{equation} \label{=t-bracket}[X,Y]_t = t[X,Y].\end{equation} 
The
function
$q(X)$ must then be changed to $q_t(X) = q(tX)$, and accordingly $\eta$
to $\eta_t$. Let
$\gu$ and
$\gv$ be compatible germs on
$\gg$, $u$ and $v$ distributions with small compact support
representing $\gu$ and $\gv$, and
$\phi$ a smooth function; we can define
\begin{equation}
\Psi(t) = \Psi_{u,v,\phi}(t) = \langle \eta_t \inv(\eta_t(u)
\ast_{G_t}
\eta_t(v)), \phi \rangle \end{equation} 
as a function of $t\in \R$. Clearly
\[\Psi(0) = \langle u\ast_{\gg} v,\phi \rangle \quad{\rm and}
\quad
\Psi(1) =\langle \eta\inv(\eta(u)\ast_G\eta(v)),\phi \rangle.\] 
The Kashiwara-Vergne conjecture is implied by (and in fact
equivalent to) the statement that for $\gu$ and $\gv$ {\em invariant}
germs, and for all $\phi$ with sufficiently small support (depending on
$u$ and $v$), the function
$\Psi_{u,v, \phi}$ is constant. 

Using the Campbell-Hausdorff formula, it can be verified that $\Psi$ is
a differentiable function of
$t$. Thus it suffices to show that for $u$ and $v$ invariant 
\begin{equation} \Psi'(t) = 0\,\, {\rm for\, all }\, \,t.\end
{equation}
In their paper
\cite{kash-ver} Kashiwara and Vergne formulate a combinatorial
conjecture on the Campbell-Hausdorff formula which implies the
vanishing of $\Psi'(t)$. It is this combinatorial conjecture which is
proven in various special cases in the papers mentioned above
(\cite{kash-ver},
\cite{rouv-81}, \cite{vergne}).
However it remains unproven in general.

Our approach is diferent. We first show that if $\gu$ and $\gv$
are invariant, then $\Psi(t)$ is {\em analytic} in $t$. Thus it
suffices to prove that  
\begin{equation} \Psi^{(n)}(0) = 0\,\, {\rm for\, all }\, \,n.\end
{equation}
While at first this may not seem to be an easier problem, however in
this paper we relate the group convolution to an extension
of the Kontsevich
$\star$-product to distributions and prove an equivalent statement
concerning the $\star$-product. 
\par
We now recall the
construction of Kontsevich. In
\cite{konts}, an associative $\star$-product is defined on any
Poisson manifold, given by a formal series in a parameter $\h$
\begin{equation} u \star_{\h} v = \sum \frac{\h^n}{n!}
B_n(u,v),
\end{equation}
where $u,v$ are $C^{\infty}$-functions on the manifold, and $B_n(u,v)$
are certain bi-differential operators. 

Consider $\gg^*$, the dual of $\gg$, equipped
with its natural Poisson structure. It is easy
to see that when $u,v$ are in $\CS(\gg)$, {\em i.e.} polynomial
functions on
$ \gg^*$, the formula for
$u\star_{\h} v$ is locally finite, so that one can set $\h =1$, and
then $u
\star_1 v$ is again in $\CS(\gg)$. Now regarding
$u$ and $v$ as distributions supported at $0$ on
$\gg$, $\star_1$ can be considered as a new convolution on
$\gg$, but defined {\em only} for distributions with point support
at $0$. 
\par
The $\star$-product is closely related to the multiplication in the
universal enveloping algebra $\CU(\gg)$, which in turn is simply
the group convolution $\ast_G$ for distributions supported at $1$ in
$G$.  Indeed, by the universal property of
$\CU(\gg)$, there exists an isomorphism between $(\CS(\gg),
\star_1)$ and $(\CU(\gg), \ast_G)$. Kontsevich shows that this
isomorphism is given explicitly in the form 
\begin{equation} u \in \CS(\gg) \mapsto \eta(u \tau\inv). \end{equation}
Here $\tau$ was defined in \cite{konts} as a formal power series
$S_1(X)$, but was shown in \cite{ads} 
to be an analytic function in a neighborhood of
$0$ in $\gg$. \footnote{A recent preprint of Shoikhet
\cite{shoi} shows that in fact $\tau
\equiv 1$. We have not used this result in our paper; its incorporation
would simplify some of the statements, but not the proofs.}
 
More generally, setting $\h$ equal to a real number $t$ we
deduce that \begin{equation} u \mapsto \eta_t(u \tau_t\inv)
\end{equation}
is an
isomorphism from $(\CS(\gg),\star_t)$ to
$(\CU(\gg),\ast_{G_t})$. This implies the identity
\begin{equation} \label{=psiconnection}
\eta_t \inv\left(\eta_t(u) \ast_{G_t}\eta_t(v)\right) =(u \tau_t \star_t
v
\tau_t) \tau_t \inv, \end{equation}
for $u,v$ distributions on $\gg$ supported at $0$.
\par
Our first main result, proved in section 3, is the following~: 

\begin{theorem}\label{konts-distr} The Kontsevich
$\star$-product on $\CS(\gg)$ extends to a (convolution) product
\begin{equation}(u,v) \mapsto u\star_t v \end{equation}
for $u$ and $v$ distributions on $\gg$ with sufficiently small support
near $0$. Moreover, formula (\ref{=psiconnection}) continues to hold.  
\end{theorem}

In section 4, we prove~: 

\begin{theorem} \label{identity-theorem} The extended $\star$-product
descends to a product on compatible germs also denoted $\star_t$. 
If 
$\gu$ and $\gv$ are compatible invariant germs, then $\gu \star_t \gv$
is invariant. Furthermore 
\begin{equation}\label{=identity} \gu \tau_t \star_t \gv \tau_t = (\gu
\ast_{\gg} \gv)\tau_t.
\end{equation}
\end{theorem}
The proof of (\ref{=identity}) requires the analyticity of $\Psi(t)$
together with
an extension of the Kontsevich homotopy argument from
\cite{ads}. 

Clearly Theorems \ref{identity-theorem} and
\ref{konts-distr} imply Theorem \ref{main-theorem}. 

{\em Acknowledgements.} Part of this research was conducted during a
visit by S.S. as {\em Professeur Invit\'e} to the
Versailles Mathematics Department (UMR CNRS 8100) and subsequently
during an NSF supported visit by
M.A. to the Rutgers Mathematics Department.

\section{Preliminaries}

Let $G$ be a finite dimensional real Lie group with
Lie algebra $(\gg,[\ ,\,])$ and fix a basis $(e_i)_{1 \leq i \leq d}$ of $\gg$.

\subsection{Symmetric Algebras}
Here, we work with the Lie algebra $\gg$ and its dual, but they are considered as vector
spaces. The symmetric algebra $\CS(\gg)$ can be considered in three different ways :
\begin{itemize}
\item as the
algebra
$\R[\gg^*]$ of polynomial functions on $\gg^*$;
\item as the algebra of
constant coefficient differential operators on $\gg$ : if $p \in \CS(\gg)$, we
write
$\partial_p$ for the corresponding differential operator. For example for $p \in \gg$,
$\partial_p$ is the constant vector field defined by $p$; 
\item finally, as the algebra (for convolution) of distributions on $\gg$ with
support
$\{0\}$ by the map $p \mapsto \bd_p$ where
\[\langle \bd_p, \phi \rangle  = \partial_p(\phi)(0) \]
for $\phi$ a test function on $\gg$. (Up to some powers of $i$, $\bd$ coincides with
the Fourier transform from functions on $\gg^*$ to distributions on $\gg$.)
\end{itemize}
When there is no ambiguity, we drop $\bd$, and use the
same notation for an element of
$S(\gg)$ as a polynomial function on $\gg^*$ and as a distribution on
$\gg$. We will then use $\cdot$ or $\ast_{\gg}$
for the product in $\CS(\gg)$ depending on how we view elements of
$\CS(\gg)$. 

The symmetric algebra $\CS(\gg^*)$ has similar interpretations.

\subsection{Algebra of differential operators} \label{s-weyl}
Let $\gW(\gg)$ the Weyl algebra of differential
operators with polynomial coefficients on $\gg$.
Any element in $\gW(\gg)$ can be uniquely written as a sum $\sum q_i
\partial_{p_i}$ with $ q_i \in \CS(\gg^*)$ and $p_i\in \CS(\gg)$. In other
words, we have a vector space isomorphism from
$\CS(\gg^*)
\otimes
\CS(\gg)$ to
$\gW(\gg)$. The inverse map associates to a differential operator in $\gW(\gg)$ its symbol
in $\CS(\gg^*)\otimes\CS(\gg)$; the symbol can be viewed as a polynomial map from $\gg^*$ to
$\CS(\gg^*)$. We observe that an element of $\gW(\gg)$ is completely determined by its
action on $\CS(\gg^*)$.
\par
Similarly, any element in $\gW(\gg^*)$
can be written as a sum $\sum_i p_i \partial_{ q_i}$ with $p_i 
\in \CS(\gg)$
and $ q_i \in \CS(\gg^*)$. 
\par
By duality with test functions, the set
of distributions
$\bD(\gg)$ on $\gg$ is a right $\gW(\gg)$-module :
\[ \langle D \cdot L, \phi \rangle  = \langle D, L \cdot \phi\rangle  \]
where $\phi \in C_c^{\infty}(\gg)$, $D \in \bD(\gg)$,
$L \in \gW(\gg)$. 
One can define
a canonical anti-isomorphism (the {\it Fourier transform})
$\CF$ from $\gW(\gg^*)$ to $\gW(\gg)$ such that
$\CF(\sum p_i \partial_{q_i}) = \sum q_i \partial_{p_i}$ for $p_i \in \CS(\gg)$ and
$q_i \in \CS(\gg^*)$. It verifies
\begin{equation} \label{=fourier on dist}\bd_{L\cdot p} = \bd_p \cdot \CF(L)
\end{equation}
for any $L \in \gW(\gg^*)$, $p\in \CS(\gg)$.
\par
\subsection{Multi-differential operators}
Let $\gW^m(\gg^*)$ be the set of $m$-differential operators with polynomial
coefficients on $\gg^*$. These are linear combinations of operators from
$C^{\infty}(\gg^*)^{\otimes m}$ to
$C^{\infty}(\gg^*)$ of the form 
\begin{equation}
(f_1 \otimes \dots \otimes f_m)\mapsto
 p \partial_{ q_1}(f_1) \dots \partial_{ q_m}(f_m) \end{equation}
where $p\in \CS(\gg)$ and $q_i \in \CS(\gg^*)$.
There is an obvious linear
isomorphism, written $A \in\CS(\gg)
\otimes
\bigotimes^m \CS(\gg^*) \mapsto \partial_A \in \gW^m(\gg^*)$. Its inverse maps a
$m$-differential operator 
$B$ to its {\em symbol} $\sigma_B$. A symbol will often be viewed as a polynomial map
from $\gg^m$ to $\CS(\gg)$ 
\par
For multidifferential operators there is no symmetry similar to the one given
by the Fourier transform $\CF$ from $\gW(\gg^*)$ to $\gW(\gg)$. Nevertheless,
in the case of bi-differential operators, for any $B = p\partial_{ q_1}
\otimes \partial_{ q_2} \in \gW^2(\gg^*)$, we define an operator $\CF(B)$
mapping functions on $\gg$ to functions on $\gg \times \gg$ :
\begin{equation}\label{co-fourier} \CF(B)(f)(x,y) =  q_1(x)
 q_2(y) [\partial_p (f)](x+y). \end{equation}
By duality, we get a ``right" action on pairs of distributions on $\gg$
in the following way~:
\begin{equation} \label{=bi-fourier} (u,v)\cdot \CF(B) =
[(u\cdot  q_1) \ast_{\gg} (v \cdot  q_2] \cdot
\partial_p ,\end{equation} where $u,v$ are distributions on $\gg$,
$\ast_{\gg}$ is the convolution on $\gg$ and it is assumed that the
convolution makes sense. 
We then have a formula similar to (\ref{=fourier on dist})
\begin{equation}\label{bi-fourier on dist}
\bd_{B(f,g)} = (\bd_f,\bd_g)\cdot\CF(B).
\end{equation}
Note that a bi-differential operator with polynomial coefficients is completely determined
by its action on point distributions. 
\par
\subsection{Enveloping algebra}
The enveloping algebra $\CU=\CU(\gg)$
of $\gg$ can be seen as the algebra of {\it left} invariant differential operators on $G$
(multiplication being composition of differential operators),
as the algebra of distributions on $G$ with support $1$,
multiplication being convolution of distributions.
Depending on the situation, we will write $\cdot$ or $\ast_G$ the product in
$\CU(\gg)$. 

It is well known, and can be easily seen, that the
symmetrization map $\beta$ from $\CS(\gg)$ to
$\CU(\gg)$ can be interpreted as the pushforward of distributions from the
Lie algebra to the Lie group by the exponential map~: for all $p\in
\CS(\gg)$, $\beta(p) = \exp_* p$. (Note that there is no Jacobian
involved here.) 
\par
\subsection{Poisson structure}
As is well known, the dual $\gg^*$ of $\gg$ is a Poisson manifold. For
$f,g$ functions on $\gg^*$, the Poisson bracket is
\begin{equation}
\{f,g\}(\nu) = \frac 12 \nu([df(\nu),dg(\nu)])
\end{equation}
for $\nu \in \gg^*$, where $df(\nu), dg(\nu) \in \gg^{**}$ are identified with
elements of $\gg$. Using the given basis $e_i$ of $\gg$, and the corresponding
structure constants
\begin{equation}
[e_i,e_j] = c_{ij}^k e_k
\end{equation}
(we use the Einstein convention of summing repeated indices).
We write $\partial_j$ for the partial derivatives with respect to the dual
basis $e_j^*$ of $e_i$. We get the following
formula for the Poisson bracket
\begin{equation}\label{=poisson}
\gamma(f,g)(\nu) = \{f,g\}(\nu) = \frac 12 c_{ij}^k \nu_k \partial_i (f)(\nu)
\partial_j (g)(\nu);
\end{equation}
said otherwise, the corresponding Poisson tensor is
$\gamma^{ij}(\nu) = \frac 12 c_{ij}^k
\nu_k$

\section{General facts about the Kontsevich \\ construction of a
$\star$-product}
\subsection{Admissible graphs}
Consider, for every integers $n\geq1$ and $m \geq 1$ a set
$G_{n,m}$ of labeled oriented graphs $\Gamma$ with
$n+m$ vertices : $1,\dots,n$ are the vertices of ``first kind", $\bar 
1,\dots,\bar
m$  are the vertices of ``second kind". Here labeled means that the edges
are labeled. We will actually use these graphs here mostly for $m = 1$ (resp.
$2$). In those cases, the vertices of second type will be named $M$ (resp.
$L,R$). The set of vertices of
$\Gamma$ is denoted $V_\Gamma$, the set of vertices of the first kind
$V^1_{\Gamma}$, and the set of edges
$E_\Gamma$. For $e = (a,b) \in E_\Gamma$, we write $a = a(e)$ and $b =
b(e)$.  For graphs in $G_{n,m}$, we assume that for any edge
$e$, $a(e)$ is a vertex of the first kind, $b(e)\neq a(e)$ and for
every vertex $a\in\{1,\dots,n\}$, the set of edges starting at $a$ has 2
elements, written $e_a^1,e_a^2$. Finally, we assume that there are no
double edges.
\par
As usual, a {\em root} of an oriented graph is a vertex which is the 
end-point of
no edge, and a {\em leaf} is a vertex which is the beginning-point of no edge.

\subsection{Multidifferential operators associated to graphs}
Let
$X$ be a real vector space of dimension
$d$ with a chosen basis
$v_1,\dots,v_d$. The vector field $\partial_{v_j}$ on $X$
is denoted $\partial_j$. We fix a $C^{\infty}$ bi-vector field $\alpha
=\sum_{i,j \in \{1,\dots,d\}} \alpha^{ij} \partial_i \partial_j$ on
$X$.
\subsubsection{Differential operators}
Let $\Gamma \in G_{n,1}$. We define a
differential operator
$D_{\Gamma,\alpha}$ by the formula~:
\begin{equation}\label{=diff}
D_{\Gamma,\alpha}(\phi) = \sum_{I\in \CL} \bigg[\prod_{k=1}^n
\big(\prod_{e \in E_{\Gamma}, b(e) = k} \partial_{I(e)}\big)
\alpha^{I(e_k^1)I(e_k^2)}\bigg] \big( \prod_{e \in E_{\gamma}, b(e) = M}
\partial_{I(e)}\big) \phi
\end{equation}
where $\CL$ is the set of maps from  $E_\Gamma$ to $\{1,\dots d\}$
(taggings of edges). One should think of this operator in the following way :
for each tagging of edges, ``put" at each vertex $\ell$ the coefficient
$\alpha^{ij}$ corresponding to the tags of the edges originating at
$\ell$ and ``put" $\phi$ at $M$. Whatever is at a given vertex
($\ell\in \{1,\dots,n\}$ or $M$) should be differentiated according to the
edges ending at that vertex.
\par
\subsubsection{Bi and multi-differential operators}
Consider a graph $\Gamma \in G_{n,2}$. Define a bi-differential
operator
$B_{\Gamma,\alpha}$ by the formula :
\begin{multline}\label{=bidiff}
B_{\Gamma,\alpha}(f,g) = \sum_{I\in \CL} \bigg[\prod_{k=1}^n
\big(\prod_{e \in E_{\Gamma}, b(e) = k} \partial_{I(e)}\big)
\alpha^{I(e_k^1)I(e_k^2)}\bigg]\\
\big( \prod_{e \in E_{\gamma}, b(e) =
L}\partial_{I(e)}\big) f \big(\prod_{e\in E_{\gamma}, b(e) = R}
\partial_{I(e)}\big) g.
\end{multline}
The bi-differential operator $B_{\Gamma,\alpha}$ has a similar 
interpretation as
$D_{\Gamma,\alpha}$ : ``put" at each vertex of the first kind the coefficient
$\alpha_{ij}$, put at $L$ the function $f$ and at $R$ the function $g$,
differentiate whatever is at some vertex according to the labels of the
edges ending at that vertex, and finally multiply everything.
\par
For graphs $\Gamma$ in $G_{n,m}$, by a straightforward generalization of
(\ref{=bidiff}), one defines $m$-differential operators, with the same
interpretation as before.
\subsection{Lie algebra case}
\subsubsection{Relevant graphs} Assume that the vector space $X = \gg^*$, and
that the bi-vector field
$\alpha$ is the Poisson bracket $\gamma$, so that
$\alpha^{ij} (\nu) = \gamma^{ij}(\nu) = \frac 12 c_{ij}^k \nu_k$. Since we
use only this bi-vector field, we shall drop $\alpha$ from the notation.
Because of the linearity of the bi-vector field associated with the Poisson
bracket, for the graphs
$\Gamma \in G_{n,m}$ for which there exists a vertex of the
first kind $\ell$ with at least two edges ending at $\ell$, $B_\Gamma = 0$.
We will call {\em relevant} the graphs which have at most one edge ending at
any vertex of the first kind; since we are dealing 
exclusively with the
Lie algebra case, we will henceforth change our notation slightly, {\em and
use the notation $G_{n,m}$ for the set of relevant graphs.}

\subsubsection{Action on distributions}
Let $\Gamma \in G_{n,2}$, and let us interpret ``graphically" the operator
$\CF(B)$ acting on pairs of distributions $u,v$ on $\gg$.  As in (\ref
{=bidiff}), the action is defined as a sum over all taggings of edges of
termes obtained in the following way : ``put"
$u$ and
$v$ at vertices $L$ and $R$, put at any vertex of the first kind the
distribution $\bd_{[e_i,e_j]}$, where  $i,j$ are the tags of the
edges originating at that vertex; any vertex tagged by $\ell$ gives a
multiplication by function $e^*_{\ell}$ of the distribution sitting at the
end point of that vertex. And finally, one takes the convolution of the
distributions at all vertices.

\subsection {The Kontsevich $\star$-product}

Kontsevich defines a certain compactification  $\o \bC_{n,m}^+$ of 
the configuration space of $n$ points $z_1,\dots,z_n$ in
the  Poincar\'e upper
half space, with 
$z_i \neq z_j$ for $i \neq j$ , and $m$ points $y_1 < y_2<\dots < y_m$ in $\R$,
up to the action of the
$az +b$-group for $a\in \R^{+*}$ and $b \in \R$. To each graph $\Gamma$ is
associated a weight
$w_{\Gamma}$ which is an integral of
somme differential form on $\o C_{n,m}^+$.

Let $X = \R^d$ with a Poisson structure $\gamma$ considered as a bi-vector
field. Let $A = C^{\infty}(X)$ the corresponding Poisson algebra.
A $\star$-product on $A$ is an associative
$\R[[\h]]$-bilinear product on $A[[\h]]$ given by a formula of the type :
\begin{equation}
\label{=star-product} (f,g) \mapsto f \star_{\h} g = fg + \h B_1(f,g) +
\frac 1{2!} \h^2 B_2(f,g) +\dots + \frac 1{n!} \h^n B_n(f,g) + \dots
\end{equation}  where $B_j$ are bi-differential operators and such that
\begin{equation}\label{=commutator}
f \star_{\h} g - g \star_{\h} f = 2\h\gamma(f,g)
\end{equation} modulo terms in $\h^2$.

The Kontsevich $\star$-product is given by formulas
\begin{equation}\label{=konts-star}
B_n(f,g) = \sum_{\Gamma \in G_{n,2}} w_{\Gamma}
B_{\Gamma}(f,g).
\end{equation}
Note that each $B_n(\cdot,\cdot)$, being a {\em finite} sum, is a
bi-differential operator.

\subsubsection{Setting
$\h=1$.} As is explained in
\cite{ads}, one can set
$\h=1$, or, for that matter, $\h = t$ for any $t\in \R$ in the Kontsevich formula
in the Lie algebra case. Indeed, for fixed
$f,g \in \CS(\gg)$ of degrees $p,q$ respectively, the terms
$B_{\Gamma,\gamma}(f,g) =0$ for $n >p+q$, so that the sum
(\ref{=star-product}) actually involves a finite number of terms. This
operation is written $\star$ rather than $\star_1$.
\par
Replace the Lie algebra $\gg$ by $\gg_{t}$ for $t \in \R$, or $t$ a formal
variable as in (\ref{=t-bracket}). Is is easy to check that the
Kontsevich $\star$-product $\star_1$ for the Lie algebra $\gg_t$ at
$\h=1$ coincides with $\star_t$.

\section{An extension of the Kontsevich
$\star$-product} 

\subsection {The symbol of the $\star$-product}

We consider here the $\star$-product as a formal bidifferential operator, 
{\em i.e.} as an element of $\gW^2(\gg^*)[[\h]]$. It has a symbol $A_{\h}$
which belongs to $\CS(\gg^*) \otimes \CS(\gg^*) \otimes \CS(\gg)[[\h]]$. 
We will view $A_{\h} = A_{\h}(X,Y)$ as a polynomial map of $(X,Y)\in \gg \times \gg$
into $\CS(\gg)$ and prove some properties of $A_{\h}$.
\par
Recall that, for $\Gamma \in G_{n,2}$, $\sigma_{\Gamma}$ is the symbol of
the corresponding bi-diffferential operator $B_{\Gamma}$.
Let
\begin{equation}\sigma_n = \sum_{\Gamma \in G_{n,2}} w_{\Gamma}
\sigma_{\Gamma}.\end{equation}
We have, for $X,Y \in \gg$
\begin{equation}A_{\h}(X,Y) = \sum \frac{\h^n}{n!}
\sigma_n(X,Y).\end{equation}
We will describe a factorization of $A_{\h}(X,Y)$ in Proposition \ref{factorization} below.
We first need to have a careful inspection of various graphs and their symbols.    

\subsubsection{Wheels}
We say that a graph $\Gamma \in G_{n,m}$ contains a {\em wheel} of length $p$ 
if there is a finite sequence
$\ell_1,\dots,\ell_p \in \{1,\dots,n\}$ with $p \geq 2$ such that
$(\ell_1,\ell_2),\dots,(\ell_{p-1},\ell_p),(\ell_p,\ell_1)$ are
edges, and there are no other edges beginning at one of the $\ell_k$ and ending
at another $\ell_{k'}$. The graph $\Gamma$ {\em is} a wheel if $p=n$. 

\subsubsection{Simple components}
Let $\Gamma \in G_{n,m}$. Consider the graph $\wt\Gamma$ whose vertices
are $\{1,\dots,n\}$ and whose edges are those edges $e\in E_{\Gamma}$ such that
$b(e) \notin {\bar 1,\dots,\bar m}$. Let $\wt\Gamma_i$ ($i \in I$) be the
connected components of $\Gamma$. The corresponding {\em simple} components
$\Gamma_i$ are the graphs whose edges are the edges of $\wt\Gamma_i$
and the edges of $\Gamma$ beginning at a vertex of $\wt\Gamma_i$, and whose
vertices are the vertices of $\wt\Gamma_i$ and $\bar 1,\dots,\bar m$.
It is easy to see that any simple component of a graph in $G_{n,m}$ can be
identified to a graph in $G_{n',m}$ for $n' \leq n$.
\par
In this situation, we use the
notation
\begin{equation}
\Gamma = \coprod_i \Gamma_i,
\end{equation}
and more generally, if $\Gamma'$ and $\Gamma''$ are sub-graphs whose simple
components determine a partition of the $\Gamma_i$, we write
\begin{equation}
\Gamma = \Gamma' \amalg \Gamma''.
\end{equation}
A {\em simple graph} is a graph with only one simple component.

\subsubsection{Symbols}
We now give rules to compute the symbol of the operator attached to a graph 
$\Gamma \in
G_{n,m}$. 
They are recorded as a series of lemmas whose easy proofs are left to the reader.

\begin{lemma} \label{degrees} Let $\Gamma \in G_{n,2}$ with $r$ roots. The symbol
$\sigma_{\Gamma}$ is of total degree $n+2r$, of  partial degree $r$ for the $\CS(\gg)$
components (we call this degree the {\em polynomial} degree) and of partial degree
$n+r$ for the $\CS(\gg^*) \otimes \CS(\gg^*)$ component (the {\em differential} degree).
\end{lemma}
\par
\begin{lemma} \label{wheel-with-leaves} Let $\Gamma \in G_{n,m}$. Assume that
there exists a sub-graph
$\Gamma_0$ with $p$ leaves $\ell_1,\dots,\ell_p$ such that
\begin{itemize}
\item each leaf $\ell_j$ is the end-point of one edge
\item $\Gamma$ is the union of $\Gamma_0$ and of $p$ sub-graphs
$\Gamma_1,\dots,\Gamma_p$
\em each $\Gamma_j$ has a single root $\ell_j$.
\end{itemize}
Then
\begin{equation}
\sigma_{\Gamma} = \sigma_{\Gamma_0}(\sigma_{\Gamma_1},\dots,\sigma_{\Gamma_p}).
\end{equation}
\end{lemma}
\par
\begin{lemma} If $\Gamma \in G_{n,n}$ is a wheel of length $n$, then
\begin{equation}
\sigma_{\Gamma}(X_1,\dots,X_n) = \tr(\ad X_1 \dots \ad X_n).
\end{equation}
\end{lemma}
\par
\begin{lemma} The symbol associated to the graph $\Gamma \in G_{1,2}$ 
whose edges
are $(1,\bar 1)$ and $(1,\bar 2)$ is the map
\[(X,Y) \in \gg \times \gg \mapsto \frac 12 [X,Y] \]
\end{lemma}
\par
\begin{lemma} Let $\Gamma \in G_{n,m}$. There exist two subgraphs
$\Gamma_1$ and $\Gamma_2$ (possibly empty), with $\Gamma_1$ having wheels and no roots, and
$\Gamma_2$ having roots and no wheels such that $\Gamma = \Gamma_1 \amalg
\Gamma_2$. The decomposition is unique up to labelling of vertices.
\label{wheel-root-decomp}
\end{lemma}
\par
\begin{lemma} Let $\Gamma \in G_{n,m}$ and assume that $\Gamma = \Gamma' \amalg
\Gamma ''$. Then
\[ \sigma_{\Gamma} =  \sigma_{\Gamma'} \sigma_{\Gamma''}. \]
\end{lemma}
\par
\subsubsection{Decomposition of $A_{\h}$}
Let us consider the following two subsets of
$G_{n,m}$ :
\begin{align*}
G_{n,2}^w &= \{\Gamma \in G_{n,2},\ \Gamma\ \rm {with\ no\ roots}\}\\
G_{n,2}^r &= \{\Gamma \in G_{n,2},\ \Gamma\ \rm {with\ no\ wheels}\},
\end{align*}
and consider the two following symbols :
\begin{align*}
A^w_{\h}(X,Y) &= \sum_n \frac{\h^n}{n!} \sum_{\Gamma \in G_{n,2}^w}
\sigma_{\Gamma}(X,Y)\\
A^r_{\h}(X,Y) &= \sum_n \frac{\h^n}{n!} \sum_{\Gamma \in G_{n,2}^r}
\sigma_{\Gamma}(X,Y).
\end{align*}

\begin{proposition} \label{factorization} 1. As an $\CS(\gg)$ valued map on $\gg
\times \gg$, $A^w_{\h}$ is scalar valued.

\noindent 2. The symbol $A_{\h}(X,Y)$ decomposes as a
product
\[A_{\h}(X,Y) = A^w_{\h}(X,Y) A^r_{\h}(X,Y).\]
\label{plustard}
\end{proposition}
\begin{proof} 1. By Lemma \ref{degrees}, the symbol of a graph with no roots has
polynomial degree $0$. This proves the first assertion.

2. By
Lemma
\ref{wheel-root-decomp}, any graph
$\Gamma$ is
$\Gamma_1 \amalg \Gamma_2$, $\Gamma_1$ with no roots, $\Gamma_2$ with no
wheels. It is easy to see from the definition that the weights are
multiplicative~:
$w_{\Gamma} = w_{\Gamma_1} w_{\Gamma_2}$. The only things that remain
to be considered are the $n!$ factors~: when one looks at a decomposition
$\Gamma = \Gamma_1 \amalg \Gamma_2$, it is unique only up to labelling, which
explains the factors $\binom n{n_1}$.
\end{proof}

For $X,Y \in \gg_{\h}$, let $Z_{\h}(X,Y)$ be their
Campbell-Hausdorff series.  Writing $Z = Z_1$, we have
\begin{equation} Z_{\h}(X,Y) = \h \inv Z(\h X, \h Y). \end{equation}
It is well known that the Campbell-Hausdorff series $Z(X,Y)$ converges for $X,Y$ small
enough. Therefore, for any
$h_0$, there exists $\epsilon$ such that for $\norm{X}$, 
$\norm{Y}\leq\epsilon$,
the power series (in $\h$) $Z_{\h}$ converges normally for $\h \leq t_0$.
\par
\smallskip
The following is due to V. Kathotia \cite[Theorem 5.0.2]{kath} :

\begin{proposition} $A^r_{\h}(X,Y) = e^{Z_{\h}(X,Y) - X - Y}$.
\label{kathotia}
\end{proposition}

As a consequence, for $\norm{X}$ and $\norm{Y}$ small enough, the formal series
$A^r_{\h}$ converges for $\h =1$.

\begin{proposition} Considered as an $\CS(\gg)$-valued function on 
$\gg \times \gg$,
$A^w_{\h}(X,Y)$ is a convergent series in a neigborhood of $(0,0)$. Moreover
$A^w_{\h}(X,Y) = \exp(A^{(w)}_{\h}(X,Y))$ where $A^{(w)}_{\h}$ is the
contribution to the series corresponding to graphs with exactly one wheel.
\end{proposition}

\begin{proof}
By an argument similar to the one in the proof of Proposition \ref{plustard},
one proves that $A_{\h}^{w}(X,Y) = \exp (A_{\h}^{(w)}(X,Y))$.

We now need to prove
convergence of the series $A_{\h}^{(w)}(X,Y)$, and by homogeneity one can
set $\h=1$. Without loss of generality, one can assume that the structure
constants $c_{ij}^k \leq 2$. Therefore
\[ \frac 1{2^p} |\tr (\ad e_{i_1} \ad e_{i_2} \dots \ad e_{i_p})| \leq d^p.\]
And more generally
\[ \frac 1{2^p} |\tr (\ad z_{1} \ad z_{2} \dots \ad z_{k})| \leq d^p.\]
with $z_{j}$ Lie monomial in $e_{i}$ of degree $p_j$ and
$\sum p_j=p$.

Let $\Gamma \in G_{n,2}^{w}$ with only one wheel of length $p \leq n$. If
the
absolute values of all components $x_i,y_j$ of $X,Y$ respectively are
less than $r$,
\begin{equation} |\sigma_{\Gamma}(X,Y)| \leq r^n d^n d^n = r^n d^{2n}.
\end{equation}

Besides, the following inequalities can be found in \cite[Lemma 2.2 and
2.3]{ads}~:
\[ |w_{\Gamma}| \leq 4^n, \ \ |G_{n,2}| \leq (8e)^n n!.\]
Finally, the terms of the series in $A_1^{(w)}(X,Y)$ can be bounded by
$(32 e)^n r^n d^{2n}$, which proves convergence for $r$ small enough.
\end{proof}
\smallskip

\subsection{A formula for the $\star$-product}
\begin{proposition} \label{konts-integral}Let $t\in \R$. Let $u,v\in \CS(\gg)$, and 
$u \star_{t} v$ their $\star$-product, considered as distributions on
$\gg$. Then $u \star_{t} v$ is given by the formula 
\begin{equation}\label{=konts-int-form}\langle u \star_t v, \phi\rangle  = 
\langle u \otimes v, A_t^w (\phi\circ Z_t)\rangle \end{equation}
for $\phi$ a test function on $\gg$.
\end{proposition}
In this proposition, we are of course looking at the Fourier transform of $\star_t$ (see
\ref{bi-fourier on dist}), meaning that we are actually expressing $\bd_{u\star_t
v}$ in terms of
$\bd_u,\bd_v$, but we avoid using these cumbersome notations. (See also \cite{arnal}
for a similar ``integral" formula.)
\begin{proof}
Since for fixed $u,v$ the series for $u
\star_{\h} v$ is finite, we can substitute to $\h$ a real
number $t$, and interchange summation with the test
function-distribution bracket~:
\begin{equation} \langle u \star_t v, \phi\rangle = \langle
u(X)\otimes v(Y), [\partial_{A_t(X,Y)} \phi](X+Y)\rangle,\end{equation}
and we avoid
convergence problems for a fixed $t$ by taking a test function
$\phi$ with small enough support.  We know that
$A_{t}(X,Y) = A_{t}^w(X,Y) A_{t}^r(X,Y)$, and
$A_{t}^w(X,Y)$ is scalar valued (Proposition
\ref{plustard}). So
$[\partial_{A_{t}(X,Y)}] \phi(X+Y) = A_{t}^w(X,Y)
[\partial_{A_{t}^r(X,Y)}\phi] (X+Y)$ By \ref{kathotia},
$A_{t}^r(X,Y) = e^{Z_{t}(X,Y) - X -Y}$, and by Taylor's formula for polynomials, 
\begin{equation}\partial_{A_{t}^r(X,Y)} \phi(X+Y) = \phi(X+Y + Z_{t}(X,Y) 
- X - Y) =
\phi(Z_t(X,Y)).\end{equation} The proposition follows.
\end{proof}

\subsection{An extension of the Kontsevich $\star$-product}

In this paragraph, we use formula (\ref{=konts-int-form}) to {\em define} the 
$\star$-product of distributions with small enough compact support, and we prove that
the definition agrees with the original one for distributions with point support. 

\begin{proposition-definition} \label{konts-gen} Let $t$ be a fixed real number. 
Then for
all distributions $u$ and
$v$ on $\gg$ with sufficiently small support, the formula (\ref{=konts-int-form}):
\[\phi \mapsto \langle u \star_t v, \phi\rangle  = 
\langle u \otimes v, A_t^w(X,Y) \phi(Z_t(X,Y))\rangle \]
defines a distribution $u \star_t v$ on $\gg$. 
\end{proposition-definition}
\begin{proof}
Let $K_u$ and $K_v$ be the supports of $u$ and $v$. Assume that $K_u$ and $K_v$
are included in a ball or radius $\alpha$. Assume that
$2 \alpha$ is less than the radius of convergence of $A_t(X,Y)$. It is clear that
(\ref{=konts-int-form}) makes sense. We need to prove that the functional $t$ defined
therein is a distribution. Since $(X,Y) \mapsto Z_t(X,Y)$ is analytic from $\gg
\times \gg$ to $\gg$, the pushforward of the compactly supported distribution $u\otimes
v$ under $Z_t$ is a distribution. But $u \star_t v$ is obtained by multiplying this
distribution by the analytic function $A^t_w$, so it is a distribution.\end{proof}

By Proposition \ref{konts-integral}, we see that Definition \ref{konts-gen} 
extends the $\star$-product. This proves
the first statement of Theorem \ref{konts-distr}.

\subsection{A connection between the $\star$-product and group convolution}

Two elements of $\CS((\gg^*)) = \R[[\gg]]$ 
play a crucial
role in this situation~:
\begin{align*}
q(X) &= \left(\det{}_{\gg} \frac{e^{\ad X/2} - e^{- \ad X/2}}{\ad X}\right)^{1/2} \\
\tau(X) &= \exp\left(\sum_n \frac {w_n} {2^n} \tr ((\ad X)^n)\right)
\end{align*}
where $w_n$ is the weight corresponding to the graph with one wheel of order
$n$. It is clear by definition that $q(X)$ is analytic on $\gg$, and it is proven in
\cite{ads} that $\tau(X)$ is analytic in a neighborhood of $0$ (See also the footnote in
the introduction). 
\par
Consider the ``infinite order constant coefficient differential operator"
$\partial_{\tau}$ with symbol $\tau$ (note that $\partial_{\tau}$ is denoted
$I_T$ in
\cite{konts}). Similarly,
$\partial_q$ (written $I_{\rm{st}}$ in \cite{konts}) is well defined.
\par
As in \cite{konts}, we write $I_{\rm{alg}}$ for the isomorphism from 
$(\CU(\gg),
\ast_G)$ to $(\CS(\gg),\star)$ coming from (\ref{=commutator}) and 
the universal
property of
$\CU(\gg)$. Recall that $\beta $ is the symmetrization map ($I_{\rm{PBW}}$ in
\cite {konts}). By \cite{konts} these four maps are related by the following
\begin{equation} I_{\rm{alg}} \inv \circ \partial_{\tau} = \beta  \circ 
\partial_q.
\label{=konts-circ}\end{equation}
When we consider elements in $\CS(\gg)$ and $\CU(\gg)$ as distributions,
$\partial_{\tau}$ and $\partial_q$ should be replaced by 
multiplication operators,
so that (\ref{=konts-circ}) is equivalent to
\begin{equation} I_{\rm{alg}}\inv (p) = \beta (p\, q \tau \inv),
\label{konts-circ-fou}
\end{equation}
and as mentioned before $\beta $ is interpreted as the direct image of distributions 
under the exponential map. We therefore have the identity (\ref{=psiconnection})~:
\[
\eta_t \inv\left(\eta_t(u) \ast_{G_t}\eta_t(v)\right) =(u \tau_t \star_t
v\tau_t) \tau_t \inv, \]
for $u,v$ distributions on $\gg$ supported at $0$. 

The following Proposition will finish the proof
of Theorem \ref{konts-distr}.

\begin{proposition} 1. The scalar valued function $A^w$ is given by
\begin{equation}\label{=awform} A^w(X,Y) = \frac{r(X) r(Y)}{
r(Z(X,Y))}
{\text \;where\;} r = q \tau \inv.\end{equation}
\noindent 2. Let $u,v $ be distributions on $\gg$ with (small
enough) compact support. The identity (\ref{=psiconnection})~:
\[ 
\eta_t \inv\left(\eta_t(u) \ast_{G_t}\eta_t(v)\right) =(u \tau_t \star_t
v\tau_t) \tau_t \inv \]
holds.\end{proposition}
\begin{proof} 
\noindent 1. Let $u,v \in \CS(\gg)$. The convolution $\beta (u r)
\ast_{G} \beta (v r)$ is given by~:
\begin{equation}\label{=group-conv}
\langle \beta (ur) \ast_{G} \beta (v r), \psi \rangle = 
\langle u r \otimes v r, \psi \circ \exp  Z\rangle, 
\end{equation}
for $\psi$ a test function on $G$. Since $I_{\rm{alg}}$ is an algebra homomorphism, by
(\ref{=konts-circ}) we get 
\begin{equation} \beta ((u\star v) r) = \beta (ur)
\ast_{G} \beta (v r).\end{equation}
Therefore, applying (\ref{=konts-int-form}), we get 
for any distributions $u,v$ with support $0$~: 
\begin{equation} 
\langle u\otimes v,  A^w (r \circ Z) \psi \circ \exp Z \rangle
	= \langle
u r \otimes v r, \psi \circ \exp Z
\rangle.
\end{equation}
Fix a neighborhood of $0$ in $\gg \times \gg$ on which $A_w(X,Y)$ is defined. Assume
that $\psi(\exp(Z(X,Y))) \equiv 1$ on that neighborhood. We have
\begin{equation} 
\langle u \otimes v, A^w (r\circ Z) \rangle
= \langle u \otimes v, r \otimes r \rangle.
\end{equation}
Since this equality holds for any $u$ and $v$ supported at $0$, this implies that the
two functions
$A^w(X,Y) r(Z(X,Y))$ and
$r(X) r(Y) $ have same derivatives at $(0,0)$; since they are analytic,
they must be equal.

\noindent 2. The second statement follows immediately from (\ref{=group-conv}), 
(\ref{=awform}) and
(\ref{=konts-int-form}).
\end{proof}

We will prove now one last property about the extended $\star$-product, relating
it directly to graphs~: 
\begin{proposition}\label{taylor} Let $u,v$ be two
distributions with small enough compact support on
$\gg$. The following holds : 
\begin{equation} {\frac{d^n}{dt^n}}\bigg|_{t=0} u \star_t v = 
\sum_{\Gamma
\in G_{n,2}} (u,v) \CF(B_{\Gamma}).\end{equation}
\end{proposition}
 
Before giving the proof, let us observe that we fall short
from proving that $u \star_t v$ is {\em analytic} in $t$ for
$u,v$ general. However, we will prove such a result for
$u,v$ invariant in section 4.   

\begin{proof}
For $u,v$ be
two distributions on
$\gg$ with small support, and $\phi$ a test function on $\gg$, let
us define
\begin{equation} \langle C_n(u,v),\phi\rangle =
\frac{d^n}{dt^n}\bigg|_{t=0} \langle u \star_t
v,\phi\rangle.\end{equation}
By inspection, $C_n$ is a bi-differential operator with
polynomial coefficient acting on distributions. 
\par
For distributions $u,v \in \CS(\gg)$, we {\em know} that
\begin{equation} C_n(u,v) = \sum_{\Gamma \in G_{n,2}} (u,v)\cdot
\CF(B_{\Gamma}).\end{equation}
Now bi-differential operators with polynomial coefficients are
completely determined by their action on distributions with
point support. This proves that 
\begin{equation} C_n =  \sum_{\Gamma \in G_{n,2}} 
\CF(B_{\Gamma}). \end{equation} 
\end{proof}

\section{Proof of Theorem \ref{identity-theorem}}

\subsection{Convolution on the level of germs}
The first step is to transfer the results of Section 3 to germs. We will
begin by giving a few definitions.
 
\subsubsection{Germs}
Recall that the {\it germ at $0$} (resp.
at $1$) of a distribution $u$ on $\gg$ (resp. $G$) is the equivalence
class of $u$ for the equivalence relation $u_1 \sim u_2$ if and
only if there exists a neighborhood $C$ of $0$ (resp.
$1$) such that for any test function
$\phi$ on $\gg$ (resp. $G$) with support in $C$,
$\langle u_1,\phi\rangle  = \langle u_2,\phi\rangle $. Clearly, for any 
distribution $u$ and any given neighborhhod $\Omega$ of $0$ (resp. $1$), there
exists a distribution with support in $\Omega$ defining the
same germ at $0$ as $u$.

\subsubsection{Action of $G$ on germs}
From the action of $G$ on $G$ by conjugation (resp. the adjoint action on 
$\gg$) we get an action of $G$ on functions and distributions on $G$ (resp.
$\gg$). For $g \in G$ and $u$ a distribution, we write 
$u^g$ for the image of $u$ under $g$. We get therefore an infinitesimal action
of $\gg$ on functions and distributions on $\gg$, which is exactly
the action by adjoint vector fields $\adj_X$ for $X \in \gg$. It is
straightforward to see that these actions go down to germs. 

\subsubsection{Invariant germs}
Invariant germs are defined as germs $\gu$ such that $\gu \cdot A
=0$ for all $A \in \gg$. By taking a basis of $\gg$, it is easy to see 

\begin{lemma} A germ $\gu$ is invariant if and only if, for any
distribution $u$ representing $\gu$, there exists an open
neighborhood $\Omega$ of $0$ such that, for any $\phi$ supported in $\Omega$
and any $A \in \gg$
\begin{equation} \langle u \cdot A, \phi \rangle = 0. \end{equation}
\end{lemma}
\subsubsection{The compatibility condition}
For $\gu$ a germ on
$\gg$, recall we have defined $C_0[\gu]$ as
the cone $C_0(\supp u)$ for any $u$ representing $\gu$. Two germs $\gu,\gv$
are compatible if
\begin{equation} C_0[\gu] \cap - C_0[\gv] = \{0\}. \end{equation}
\subsubsection{Proof of the first statement of Theorem
\ref{identity-theorem}} We need to prove that the $\star$-product
descends to germs. The analogous statements about the convolutions
on $\gg$ and on $G$  are made in \cite{kash-ver}. By formula
(\ref{=psiconnection}), we deduce it for $u \star v$.

\subsection{Invariant germs}

In order to finish the proof of Theorem \ref{identity-theorem}, we
need to have very precise statements about the choice of
representatives of invariant germs that we will work with.  This
is what we do in this paragraph.  

We choose some norm $\Vert \, \Vert$ on $\gg$, and write $B(0,r)$ for the open
ball of center $0$ and radius $r$. For us, open cones mean cones
containing $0$ with open cross-section.

\begin{lemma} \label{choice} 1. Let $\gu$ be a germ. For any open
cone
$D$ containing $C_0[\gu]$, there exists a representative of $\gu$
with support in $D$.

2. Let $\gu,\gv$ be compatible germs. There exist open cones $D_0[\gu] \supset C_0[\gu]$,
$D_0[\gv] \supset C_0[\gv]$ such that 
\begin{equation}
\o{D_0[\gu]} \cap - \o{D_0[\gv]} = \{0\}.\end{equation}
\end{lemma}

\begin{proof} The second statement is easy. We prove the first. 
Let $u_1$ be any representative of $\gu$. By definition 
of $C_0[\gu]$, there exists $\eta >0$
such that 
\begin{math}\supp u_1 \cap B(0,\eta) \subset D.
\end{math}
Let $\chi$ be a $C^{\infty}$ function supported in $B(0,\eta)$ which is
identically equal to $1$ in $B(0,\eta/2)$. Clearly, $u = u_1 \chi$ is
a representative of $\gu$ with support included in $D_0[\gu]$.
\end{proof}
For $\beta >0$, we write $D_0^{\beta}[\gu] = D_0[\gu] \cap B(0,\beta)$. The following
lemma is crucial for our purposes.
\begin{lemma}\label{campbell-estimate} Let $\gu,\gv$ be compatible germs, and $D_0[\gu], D_0[\gv]$ open 
cones as in Lemma \ref{choice}. There exists a $\beta>0$ such that, for any $\gamma >0$,
there exists
$\delta >0$ satisfying
\begin{equation}\label{=in1}
 \supp(\phi \circ Z_t)
\cap  (D_0^{\beta}[\gu] \times 
D_0^{\beta}[\gv])\subset B(0,\gamma) \times
B(0,\gamma)\end{equation}
for all smooth $\phi$ with support in $B(0,\delta)$.
\end{lemma}
\begin{proof}
We study the restriction of $Z_t(X,Y)$ to $D_0[\gu] \cap D_0[\gv]$. The
Campbell-Hausdorff formula implies that there exists a positive number
$\beta_1$ such that the $\gg$-valued map $(t,X,Y) \mapsto Z_t(X,Y)$ is
analytic for $t \leq 1, \Vert X \Vert \leq \beta_1, \Vert Y \Vert \leq
\beta_1$. Furthermore
\begin{equation} \frac{\partial Z_t}{\partial X}(0,0) =
\frac{\partial Z_t}{\partial Y}(0,0) = I. \end{equation}
Therefore, the implicit equation $Z_t(X,Y)=0$ can be solved : there
exists a constant $\beta_2 < \beta_1$ and an analytic map $(X,t) \mapsto
z_t(X)$, defined for
$t \in [-1,1]$ and $\Vert X \Vert \leq \beta_2$ such that
$Z_t(X,Y) = 0$ is equivalent to
$Y = z_t(X)$ for $\Vert X \Vert \leq \beta_2, \Vert Y \Vert \leq \beta_2$.
Now since $\partial z_t/\partial X(0) = -I$ and 
$D[\gu]$ is an open cone, we conclude that there exists
$\beta < \beta_2$ such that, for $X \in D_0^{\beta}[\gu]$, and any $t \in
[-1,1]$,
$z_t(X)\in - D_0[\gu]$. Finally, we have proven that, for
$X,Y$ with norms less than $\beta$, $X \in D_0[\gu], Y \in D_0[\gv]$,
$Z_t(X,Y) =0$ implies $X = Y = 0$. 
Using a straightforward topological argument, we get
\begin{equation}\begin{aligned}\label{=in0}
\exists \beta, \forall \gamma \in (0,\beta], &\exists
\delta >0, \forall X \in D_0^{\beta}[\gu], \forall Y
\in D_0^{\beta}[\gv],\\
& \forall t \in [-1,1], \Vert Z_t(X,Y)\Vert \leq\delta 
{\, \rm implies \,} \Vert X \Vert \leq \gamma,\Vert Y \Vert \leq
\gamma.\end{aligned}\end{equation}
Assume that $\phi$ is a test function supported in $B(0,\delta)$. Then,
by (\ref{=in0}), 
\[ \supp(\phi \circ Z_t) \cap 
(D_0^{\beta}[\gu] \times  
D_0^{\beta}[\gv])\]
is included in $B(0,\gamma) \times B(0,\gamma)$, as we wanted. 
\end{proof}
We can now make the following precise statement
about the choice of representatives of germs.
\begin{proposition}\label{germtodist}
Let $\gu,\gv$ be compatible germs, and $D_0[\gu], D_0[\gv]$ open 
cones as in Lemma \ref{choice}.
There exists a positive real $\beta$ such that, for any $t \in [-1,1]$,
the germ at $0$ of $u \star_t v$ does not depend on the
choice of $u$ (resp. $v$) distribution supported in 
$D_0^{\beta}[\gu]$ (resp. $D_0^{\beta}[\gv]$) representing
$\gu$ (resp. $\gv$). This implies that the germ of $u \star v$ is independent
of the choice of $D_0[\gu], D_0[\gv]$.
\end{proposition}
\begin{proof} Let $\gu, \gv$ be compatible
germs, and $(u_1,v_1,),(u_2,v_2)$ two pairs of representatives. By formula 
(\ref{=konts-int-form}),
since the multiplication factor $A_t^w$ does not play a role, it is enough to
prove that, for any
$t \in [-1,1]$, the formula
\begin{equation} \label{=circ}\langle u_1 \otimes v_1,\phi \circ Z_t\rangle = 
\langle u_2 \otimes v_2,\phi\circ Z_t\rangle\end{equation}
provided the supports of $u_i,v_i,\phi$ are adequate. 
Assume that 
$(u_1,u_2)$ (resp. $(v_1,v_2)$) are
supported in
$D_0^{\beta}[\gu]$ (resp. $D_0^{\beta}[\gv]$).
Since $u_1 \sim u_2$, $v_1 \sim v_2$, there exists $\gamma$ such that for
any test function $\psi$ with support included in the ball $B(0,\gamma)$,
\begin{equation}\begin{aligned}\label{=in2} \langle u_1, \psi \rangle =
\langle u_2,\psi\rangle\\ 
\langle v_1, \psi \rangle = \langle v_2, \psi \rangle.
\end{aligned}\end{equation}
Using (\ref{=in1}) we deduce (\ref{=circ}) for $\phi$
supported in $B(0,\delta)$.
\end{proof}

\begin{proposition}\label{invariance} Let $\gu, \gv$ be two
compatible invariant germs. The germ $\gu \star_t \gv$ is invariant.
\end{proposition}
\begin{proof}
It is enough to do it for $t=1$. We consider $\beta$ from Proposition
\ref{germtodist} As before, we chose representatives
$u,v$ of
$\gu,\gv$ supported in $D_0^{\beta}[\gu]$ and $D_0^{\beta}[\gv]$ 
respectively. Since the
germs are invariant, there exists a 
$\gamma < \beta $ such that, for any test function $\psi$ supported in
$B(0,\gamma)$ 
\begin{equation}\label{=invariance}
\langle u \cdot A, \psi \rangle = \langle v \cdot A, \psi \rangle = 0
\end{equation}
for all $A \in \gg$.
Applying again Lemma \ref{campbell-estimate}, consider $\phi$ a test function
supported in $B(0,\delta)$. We prove that
\[
\langle (u \star v) \cdot A, \phi \rangle =  \langle (u \star v),
A \cdot \phi \rangle = 0
\]
for all $A \in \gg$.
Indeed, since the
functions $q,\tau$ are invariant, it is enough to prove
\begin{equation}
\langle u\otimes v, (A \cdot \phi)\circ Z)\rangle = 0.
\end{equation}
Using the covariance of $Z(X,Y)$ under the adjoint action of $g \in G$: 
\[ g.(Z(X,Y)) = Z(g.X,g.Y)\]
writing $g = \exp(tA)$ and differentiating at $t=0$, we get
\begin{equation} 
(A \cdot \phi) \circ Z  = (A \otimes 1 + 1 \otimes A) (\phi \circ Z).
\end{equation}
Thus we get
\[
\langle u\otimes v, (A \cdot \phi)\circ Z)\rangle = \langle u \cdot A \otimes v,
\phi\circ Z\rangle +\langle u \otimes v \cdot A,
\phi\circ Z\rangle.
\]
By Lemma \ref{campbell-estimate}
\[\supp(\phi \circ Z) \cap 
(D_0^{\beta}[\gu] \times D_0^{\beta}[\gv]) \subset
B(0,\gamma) \times B(0,\gamma).\]
Now we use (\ref{=invariance}) to conclude.
\end{proof}

\subsection{End of proof of Theorem \ref{identity-theorem}}
As before, we consider two compatible invariant germs $\gu$ and $\gv$, and take
representatives $u,v$ of $\gu,\gv$.  We shall
prove the equivalence of distributions~:
\begin{equation} u \tau_t \star_t v \tau_t \sim (u \ast_{\gg} v) \tau_t,  
\end{equation}
 for $u,$ and $v$
adequately chosen. It clearly will imply Theorem \ref{identity-theorem}.
It is enough to prove that for any arbitrary test
function $\phi$ with sufficiently small support
\begin{equation}\label{=final} \langle u \tau_t \star_t v \tau_t, \phi \rangle =
\langle (u \ast_{\gg} v)\tau_t, \phi. \rangle
\end{equation}
Considering both sides as functions of $t$, we will prove that they are
analytic in $t$, and that their derivatives at $0$ at any order are equal. 

\subsubsection{Analyticity}
We will derive the required analyticity result by extending some arguments of
Rouvi\`ere from \cite{rouv-90, rouv-91}. For the reader's convenience we summarize
below the results that we need from these papers.
\par
Rouvi\`ere defines an analytic function $e(X,Y)$ on $\gg \times \gg$ and a family
of maps $\Phi_t$ (depending smoothly on $t \in [0,1]$) from $\gg \times \gg$ to $\gg
\times
\gg$ with the following properties
\begin{itemize}
\item $\Phi_0 = I$.
\item For all $t$, $\Phi_t(0,0) = (0,0)$ and $\Phi_t$ is a local diffeomorphism at
$(0,0)$.
\item If $\sigma : \gg \times \gg \to \gg$ is the addition map, then $\sigma \circ
\Phi_t \inv = Z_t$.
\end{itemize} 
For a smooth function $g$ on $\gg \times \gg$, define
\begin{align*}
g_t &= (f_t g) \circ \Phi_t \inv {\rm \; where}\\
f_t(X,Y) &= \frac {q(tX) q(tY)}{q(tX + tY)} e(tX,tY) \inv.
\end{align*}
Rouvi\`ere proves
\begin{equation} \label {=r-equality}
\frac{\partial}{\partial t} g_t  = \tr{}_{\gg}\left( \ad X \circ 
\frac{\partial}{\partial X} (g_t F_t) + \ad Y \circ \frac{\partial}{\partial Y} 
(g_t G_t)\right)
\end{equation}
where $F_t$ and $G_t$ are certain smooth  functions on $\gg \times \gg$, and all
differentials are taken at $(X,Y)$. Since the differential operators on the right
can be expressed in terms of adjoint vector fields, Rouvi\`ere uses this to conclude
that, for invariant distributions
$u$ and $v$, 
\begin{equation}\label{=whence} \langle u \otimes v,  \frac{\partial}{\partial t}
g_t
\rangle
\equiv 0 {\rm \, \, whence\, \,} \langle u \otimes v, g_0 \rangle = \langle
u \otimes v, g_1 \rangle.
\end{equation}
If now $\phi$ is a smooth function on $\gg$, applying this result to 
\[g(X,Y) = \frac{e(X,Y)} {q(X)q(Y)} \phi(X+Y),\]
and pairing with $u q$ and $vq$, we obtain 
\begin{equation} \langle u \otimes v, e (\phi \circ \sigma) \rangle
= \langle u q\otimes v q,  (\phi q\inv)\circ Z \rangle.
\end{equation}
Rewriting this, we get
\begin{equation}\label{=rouviere}
 \langle u \otimes v, e (\phi \circ \sigma) \rangle =\langle \eta \inv (\eta(u)
\ast_G\eta(v)),\phi\rangle.
\end{equation}

\begin{lemma} Let $\gu,\gv$ be two compatible invariant germs, $D_0[\gu],D_0[\gv]$
chosen as in Lemma \ref{choice}, and 
$\beta$ from (\ref{=in1}). Let   
$u,v$ be  representatives
which verify the conditions of Proposition \ref{germtodist}.
There exists a positive number $\delta$ such that, for
any test function $\phi$ supported in $B(0,\delta)$, 
\begin{equation} \langle u \tau_t \star_t v \tau_t, \phi \rangle \end{equation}
is an analytic function of $t$ in a neighborhood of $[0,1]$.
\end{lemma}
Note that this lemma implies that the function $\Psi_{u,v\phi}(t)$ considered in the
introduction is analytic.
\begin{proof} The first step is to prove that (\ref{=rouviere}) holds under the
assumptions of the lemma, provided $\phi$ has small enough support. 
We proceed as in Proposition \ref{invariance}~: from the invariance of
$\gu,\gv$, we have a constant $\gamma$ such that
\begin{equation}\label{=only} \langle u \cdot A, \psi \rangle = \langle v \cdot A,
\psi
\rangle = 0 \end{equation}
for all $A \in \gg$ and any $\psi$ supported in $B(0,\gamma)$. We derive
$\delta$ from Lemma
\ref{campbell-estimate}. Let
$\phi$ have support in
$B(0,\delta)$.
In Rouvi\`ere's proof which was just outlined above, 
the distribution $u \otimes v$ is paired with the function $g_t$ and some
derivatives of $g_t$. Since $\sigma \circ \Phi_t
\inv = Z_t$, we can write
\[g_t =   K_t (\phi \circ Z_t) ,\] 
where $K_t$ is some smooth function. Therefore the support of $g_t$ (and
any derivative) is included in the support of $\phi \circ Z_t$. By (\ref{=in1}) 
wee see that 
\begin{equation} \supp g_t \cap (
D_0^{\beta}[\gu] \times 
D_0^{\beta}[\gv])\end{equation}
is included in $B(0,\gamma) \times B(0,\gamma)$. 
Since we are pairing with $u \otimes v$, with $u,v$ verifying (\ref{=only}), we still
conclude as in Rouvi\`ere that $\langle u \otimes v, g_t \rangle$ is independent of $t$, so
that (\ref{=whence}) and (\ref{=rouviere}) still hold.  

Now we apply (\ref{=rouviere}), but to the algebra $\gg_t$ (note that this $t$
is not ``the same" as the $t$ used before !). Using (\ref{=psiconnection}), since
the $e$ function for $\gg_t$ is given by $e_t(X,Y) = e(tX,tY)$ (see
\cite[Proposition 3.14]{rouv-90}) we derive~:
\begin{equation}\label{=rouv-t} 
\langle u \tau_t \star_t v \tau_t, \phi \rangle = \langle
u(X) \otimes v(Y), e(tX,tY) \tau(tX+tY) \phi(X+Y) \rangle.
\end{equation} 

Using the fact that $e, \tau, q$ are analytic in a neighborhood of $0$, we conclude
from (\ref{=rouv-t}) that $\langle u \star_t v, \phi \rangle$ is analytic for $t
\in [-1,1]$.
\end{proof}

The analyticity of the right hand side of (\ref{=final}) being 
immediate, it now remains to prove the equality of derivatives of both sides of 
(\ref{=final}) at
$0$ to all orders.

\subsubsection{Cancellation}

\begin{lemma} There exists a bidifferential operator $M_n \in \gW^2(\gg^*)$ such that
\begin{equation} (u,v)\cdot \CF(M_n) = \frac {d^n}{dt^n}\bigg|_{t=0} u \tau_t
\star_t v \tau_t - (u\ast_{\gg} v)\tau_t \end{equation}
for $u,v$ distributions with compact support. \end{lemma}
\begin{proof}
Let $\tau_{(p)}$ be the differential of $\tau$ at $0$ of order $p$,
evaluated at the $p$-tuple $(X,\dots,X)$.
It follows from Proposition \ref{taylor} that the right hand side of the lemma is
\begin{equation}\label{=connection}
\sum_{p+q+r = n} \frac
{n!}{p!\,q!\,r!} (u \tau_{(p)}, v \tau_{(q)}) \cdot B_r - (u \ast_{\gg} v)
\tau_{(n)}
\end{equation}
where
\[ B_r = \sum_{\Gamma \in G_{r,2}} B_{\Gamma}. \]
Clearly, this expression is
the Fourier transform of a bidifferential operator with polynomial coefficients.
\end{proof}

We now define a certain subset $\gJ$ of $\gW^2(\gg^*)$.
Considering the natural surjective map $\CC$ from
$\gW(\gg^*)
\otimes
\gW(\gg^*)$ to
$\gW^2(\gg^*)$ given by~:
\[p_1 \otimes q_1 \otimes p_2 \otimes q_2 \mapsto p_1p_2 \otimes q_1 \otimes q_2\]
for $p_i \in \CS(\gg)$, $q_i \in \CS(\gg^*)$. Let $B = \CC(D_1\otimes D_2)$. Then
\begin{equation}  \langle (u,v) \cdot \CF(B), \phi \rangle = \langle u \cdot \CF(D_1)
\otimes v\cdot \CF(D_2), \phi \circ \sigma \rangle.  \end{equation} Let now $\gI$ be the
left ideal in $\gW(\gg^*)$ generated by adjoint vector fields; we define 
\begin{equation}\gJ = \CC(\gI
\otimes\gW(\gg^*) + \gW(\gg^*) \otimes \gI).\end{equation} 

We now establish

\begin{lemma} For any $B \in \gJ$, and $u,v,\phi$ as in Proposition
\ref{invariance},
\begin{equation} \langle(u,v) \CF(B), \phi \rangle = 0. \end{equation}
\end{lemma}
\begin{proof}
By symmetry, it
is enough to prove that for $B = \CC(D_1 A \otimes D_2)$ with $A$ an adjoint vector field,
$D_1 \in\gW(\gg^*), D_2\in \gW(\gg^*)$, 
\begin{equation} \langle (u,v) \cdot \CF(B), \phi \rangle = 0.  \end{equation} 
But 
\begin{equation}  \langle (u,v) \cdot\CF(B), \phi \rangle = \langle u \cdot A
\otimes v, \CC(D_1 \otimes D_2) (\phi \circ \sigma)\rangle.\end{equation}
Since the support $\CC(D_1 \otimes D_2) (\phi \circ \sigma)$ is included in the support
of $\phi \circ \sigma$, we can now apply Proposition \ref{invariance}.
\end{proof}

To conclude the proof of Theorem \ref{identity-theorem}, 
we need to prove that $M_n$ belongs to $\gJ$. We adapt the proof of Theorem 0.3 of 
\cite{ads}, which relies on an argument of homotopy. We
use the notations of \cite{ads}. 

As a consequence of Stokes' formula, one can express  
\begin{equation}(u,v) \mapsto (u,v)\cdot \CF(M_n) \end{equation}
as a weighted sum of bi-differential operators
$B_{\Gamma}$ for $\Gamma \in
G_{n,2}$ with weights
\begin{equation} w'_{\Gamma} = \int_{Z_n} \omega'_{\Gamma}, \end{equation}
where $Z_n$ is some subset of the boundary of the compactification $\o
C_{n+2,0}^+$, and $\omega'_{\Gamma}$ is a differential form. Among the
configurations 1 to 4 of
\cite[Proof of Theorem 0.3]{ads} we only need to consider the two following 
configurations, the other cases being adressed in exactly the same way as in \cite{ads}~:
\begin{itemize}
\item Two-point clusters of type 1 correspond to a bi-differential operator
of the form $\CC(\gI \otimes \gW(\gg^*))$
\item Two point clusters of type 2 1 correspond to a bi-differential operator
of the form $\CC(\gW(\gg^*)\otimes \gI)$. 
\end{itemize}
\par
This finishes the proof of Theorem \ref{identity-theorem}.

\end{document}